\documentclass{article}

\usepackage{amssymb, amsmath, amsfonts}
\usepackage{mathscinet} 

\usepackage{tikz}
\usetikzlibrary{decorations,decorations.markings}
\tikzset {
  ->-/.style = {
      decoration={markings, mark=at position 0.5 with {\arrow{>}}},
      postaction={decorate}
    },
  ->--/.style = {
      decoration={markings, mark=at position 0.2 with {\arrow{>}}},
      postaction={decorate}
    },
  -<-/.style = {
      decoration={markings, mark=at position 0.5 with {\arrow{<}}},
      postaction={decorate}
    },
  x=1mm,
  y=1mm,
  baseline=-0.5ex
}

\usepackage{array}
\newcolumntype{O}{>{\centering\arraybackslash\(}p{1in}<{\)}}
\newcolumntype{T}{>{\centering\arraybackslash\(}p{2in}<{\)}}
\newcolumntype{F}{>{\centering\arraybackslash\(}p{4in}<{\)}}

\usepackage{amsthm}
\newtheorem{thm}{Theorem}
\newtheorem{prop}{Proposition}
\theoremstyle{remark}
\newtheorem*{rmk}{Remark}

\newcommand{\qbinom}[2]{\genfrac{[}{]}{0pt}{}{#1}{#2}}
\DeclareMathOperator{\sign}{sign}

\title{A new approach to the $SL_n$ spider}

\author{Stephen Bigelow}

\begin{document}

\maketitle

\begin{abstract}
  The $SL_n$ spider gives a diagrammatic way to encode the representation
  category of the quantum group $U_q(\mathfrak{sl}_n)$.  The aim of
  this paper is to define a new spider that contains the $SL_n$ spider.
  The new spider is defined by generators and relations, according to
  simple rules that start with combinatorial data coming from the
  root system of $SL_n$.
\end{abstract}

\section{Introduction}

A {\em spider} is an algebraic structure that uses diagrams to
encode the representation category of a quantum group.  Kuperberg
\cite{kuperberg} gave generators and relations for spiders for rank
two Lie algebras.  Cautis, Kamnitzer, and Morrison \cite{ckm} did
the same for $SL_n$.  The aim of this paper is to define a new
spider that contains a copy of the $SL_n$ spider.

My motivation for this construction was to improve my understanding
of spiders.  The presentations for spiders have relations that seem
far from obvious.  Our new spider has a simple presentation
defined using the combinatorics of roots and weights of $SL_n$.  I
will discuss this, as well as some hints of connections to other
approaches, in Section~\ref{sec:motivation}.

We will use a version of the ``more natural (if less formal)''
definition of spiders from \cite{kuperberg}.  For us, a \emph{web}
will be a graph drawn in the disk, with some endpoints on the
boundary of the disk.  Edges have orientations and labels, and
bivalent vertices have \emph{tags}.  Edges are also called
\emph{strands}.  Webs are taken up to boundary-fixing isotopies.

Our spiders also contain formal linear combinations of webs, with
scalars in the field $\mathbb{C}(q)$.  We reserve the word ``web''
for a single graph in the disk, as opposed to any other linear
combination.

A presentation of a spider consists of a set of labels of strands,
a set of possible vertices, and a list of linear relations between
webs.  A relation can be applied ``locally'' inside a larger web,
which remains unchanged outside the disk where the relation is
applied.

Kuperberg \cite{kuperberg} requires a web to have a distinguished
``first'' endpoint on the boundary, but we can safely ignore this.
In any relation we write between webs, all terms have the same
arrangement of endpoints on the boundary, so we could give them the
same arbitrary first endpoint.

We will give presentations of two spiders.  The first is the $SL_n$
spider from \cite{ckm}.  The webs of the second spider will be
called \emph{cobwebs}, to help distinguish them from the \emph{webs}
for the $SL_n$ spider.  We then describe an algorithm that takes a
web and outputs a sum of cobwebs.  We prove that this gives a
well-defined and injective map from the $SL_n$ spider into the
spider of cobwebs.

For an integer $n \ge 0$, the quantum integer $[n]$ is
  $$[n] = q^{n-1} + q^{n-3} + \dots + q^{3-n} + q^{1-n}.$$
For $0 \le k \le n$, the quantum factorial, and quantum binomial are
  $$[n]! = [n][n-1]\dots[1], \qquad
    \qbinom{n}{k} = \frac{[n]!}{[n-k]![k]!}.$$
We also define negative quantum integers by $[n] = -[-n]$.

\section{Webs}

\begin{figure}
\[
\begin{array}{OOOO}
  \tikz {
    \draw[->] (0,0) -- (90:6) node [right] {$k + l$};
    \draw[->-] (210:6) node [left] {$k$} -- (0,0);
    \draw[->-] (330:6) node [right] {$l$} -- (0,0);
  }
  & \tikz {
    \draw[->-] (90:6) node [right] {$k + l$} -- (0,0);
    \draw[->] (0,0) -- (210:6) node [left] {$k$};
    \draw[->] (0,0) -- (330:6) node [right] {$l$};
  }
  & \tikz {
    \draw[->] (0,1) -- (0,6) node [right] {$k$};
    \draw[->] (0,1) -- (0,-4) node [right] {$n - k$};
    \draw (0,1) -- (2,1);
  }
  & \tikz {
    \draw[->-] (0,6) node [right] {$k$} -- (0,1);
    \draw[->-] (0,-4) node [right] {$n - k$} -- (0,1);
    \draw (0,1) -- (2,1);
  }
\end{array}
\]
\caption{The four kinds of vertex in a web.}
\label{fig:webgen}
\end{figure}

In a web for the $SL_n$ spider, each strand is oriented and is
labeled by an integer from $\{1,\dots,n-1\}$.  A web can have
bivalent and trivalent vertices.  At a trivalent vertex, either
strands labeled $k$ and $l$ fuse to a strand labeled $k + l$, or a
strand labeled $k + l$ divides into strands labeled $k$ and $l$.
At a bivalent vertex, strands labeled $k$ and $n - k$ either both
enter or both exit the vertex, and a {\em tag} points to one of the
two sides of the vertex, breaking its rotational symmetry.  See
Figure~\ref{fig:webgen}.

\begin{figure}
\[
  \tikz {
    \draw[->] (0,0) -- (90:6) node [right] {$n$};
    \draw[->-] (210:6) node [left] {$k$} -- (0,0);
    \draw[->-] (330:6) node [right] {$n-k$} -- (0,0);
  }
  =
  \tikz {
    \draw (0,-1) -- (0,1);
    \draw[->-] (210:6) node [left] {$k$} ..
       controls (210:5) and (-2,-1) .. (0,-1);
    \draw[->-] (330:6) node [right] {$n-k$} ..
       controls (330:5) and (2,-1) .. (0,-1);
  }
  \qquad
  \tikz {
    \draw[-<-] (0,0) -- (90:6) node [right] {$n$};
    \draw[<-] (210:6) node [left] {$k$} -- (0,0);
    \draw[<-] (330:6) node [right] {$n-k$} -- (0,0);
  }
  =
  \tikz {
    \draw (0,-1) -- (0,1);
    \draw[<-] (210:6) node [left] {$k$} ..
       controls (210:5) and (-2,-1) .. (0,-1);
    \draw[<-] (330:6) node [right] {$n-k$} ..
       controls (330:5) and (2,-1) .. (0,-1);
  }
\]
\caption{Convention for removing a strand labeled $n$.}
\label{fig:tag}
\end{figure}

It is convenient to also allow strands labeled $0$ or $n$, with a
convention for deleting them.  Any trivalent vertex that had a
strand labeled $0$ simply ceases to be a vertex.  Any trivalent
vertex that had a strand labeled $n$ becomes a bivalent vertex,
with a tag on the side of the deleted strand.  See Figure~\ref{fig:tag}.

\begin{rmk}
  The definition in \cite{ckm} contains a minor sign error, which
  has no consequences for the rest of their results.  Their
  \cite[relations~2.3--10]{ckm} are correct.  Their convention in
  \cite[Section~2.1]{ckm} for deleting strands labeled $n$ is
  different from our Figure~\ref{fig:tag}.  With their convention,
  they cannot necessarily reverse all arrows in a relation, but
  should instead simultaneously reverse all arrows and flip all
  tags.  With our convention, we can reverse all arrows in a relation,
  without changing the tags.  Flipping the tags can either switch
  the sign of a web or have no effect.
\end{rmk}

\begin{figure}
  \[
    \begin{array}{|T|T|}
      \hline
      {
        \phantom{ \tikz {\draw (0,-10) -- (0,10)}} 
        \tikz {
          \draw[->] (0,0) -- (0,6) node [right] {$k$};
          \draw[->] (0,0) -- (0,-6) node [right] {$n - k$};
          \draw (0,0) -- (2,0);
        }
        = (-1)^{k(n - k)} \tikz {
          \draw[->] (0,0) -- (0,6) node [right] {$k$};
          \draw[->] (0,0) -- (0,-6) node [right] {$n - k$};
          \draw (0,0) -- (-2,0);
        }
      } &
      {
        \tikz {
          \draw[->-] (0,-6) -- (0,-2);
          \draw (0,-2) -- (2,-2);
          \draw[->-] (0,2) -- (0,-2);
          \draw (0,2) -- (2,2);
          \draw[->] (0,2) -- (0,6) node [right] {$k$};
        }
        = \, \tikz {
          \draw[->] (0,-6) -- (0,6) node [right] {$k$};
        }
      } \\ \hline
      {
        \phantom{ \tikz {\draw (0,-10) -- (0,10)}} 
        \tikz {
          \draw[->-] (-6,-6) node [left] {$l$}
            .. controls (-5,-5) and (-4,-3) .. (0,-3);
          \draw[->-] (6,-6) node [right] {$m$}
            .. controls (5,-5) and (4,-3) .. (0,-3);
          \draw[->-] (0,-3) -- (0,3);
          \draw[->-] (-6,6) node [left] {$k$}
            .. controls (-5,5) and (-5,3) .. (0,3);
          \draw[->] (0,3) .. controls (4,3) and (5,5) .. (6,6);
        }
        = \tikz {
          \draw[->-] (-6,-6) node [left] {$l$}
            .. controls (-5,-5) and (-3,-4) .. (-3,0);
          \draw[->-] (-6,6) node [left] {$k$}
            .. controls (-5,5) and (-3,4) .. (-3,0);
          \draw[->-] (-3,0) -- (3,0);
          \draw[->-] (6,-6) node [right] {$m$}
            .. controls (5,-5) and (3,-4) .. (3,0);
          \draw[->] (3,0) .. controls (3,4) and (5,5) .. (6,6);
        }
      } &
      {
        \tikz {
          \draw[->-] (0,-8) -- (0,-4);
          \draw[->-] (0,-4) arc (270:90:4);
          \draw[->-] (0,-4) arc (-90:90:4) node[pos=0.5,right]{$1$};
          \draw[->] (0,4) -- (0,8) node [right] {$k$};
        }
        = [k] \: \tikz {
          \draw[->] (0,-8) -- (0,8) node [right] {$k$};
        }
      } \\ \hline
    \multicolumn{2}{|F|}
    {
      \phantom{ \tikz {\draw (0,-10) -- (0,10)}} 
      \tikz {
        \draw[->-] (45:4) arc (45:135:4) node [pos=0.5, above]{$1$};
        \draw[-<-] (135:4) arc (135:225:4);
        \draw[->-] (225:4) arc (225:315:4) node [pos=0.5, below]{$1$};
        \draw[->-] (315:4) arc (315:405:4);
        \draw[->-] (7,-7) node [right] {$1$} -- (315:4);
        \draw[->] (45:4) -- (7,7) node [right] {$1$};
        \draw[->] (135:4) -- (-7,7) node [left] {$k$};
        \draw[->-] (-7,-7) node [left] {$k$} -- (225:4);
      }
      =
      \tikz {
        \draw[->-] (-7,-7) node [left] {$k$}
          .. controls (-6,-6) and (-3,-5) .. (0,-4);
        \draw[->-] (7,-7) node [right] {$1$}
          .. controls (6,-6) and (3,-5) .. (0,-4);
        \draw[->-] (0,-4) -- (0,4);
        \draw[->] (0,4) .. controls (-3,5) and (-6,6)
          .. (-7,7) node [left] {$k$};
        \draw[->] (0,4) .. controls (3,5) and (6,6)
          .. (7,7) node [right] {$1$};
      }
      + [k - 1] \tikz {
        \draw[->] (-7,-7) .. controls (-4,-4) and (-4,4)
          .. (-7,7) node [left] {$k$};
        \draw[->] (7,-7) .. controls (4,-4) and (4,4)
          .. (7,7) node [right] {$1$};
      }
    } \\ \hline
\end{array}
\]
\caption{Web relations. Some labels are omitted, but can be deduced.}
\label{fig:webrel}
\end{figure}

As they note in their second remark, \cite[relations~2.3--10]{ckm}
contains redundancies.  In this paper, we use the relations shown
in Figure~\ref{fig:webrel}, together with their mirror images and
arrow reversals.  Here, and throughout the paper, we omit some of
the labels on strands, as long as they can be deduced from context.

In the order shown in Figure~\ref{fig:webrel},
we name our relations:
  \begin{itemize}
    \item ``tag switching'' is \cite[relation~2.3]{ckm},
    \item ``tag cancelation'' is \cite[relation~2.5]{ckm}
      in the special case $k+l=n$,
    \item ``$I = H$'' is \cite[relation~2.6]{ckm},
    \item ``bursting a digon'' is \cite[relation~2.4]{ckm}
      in the special case $l=1$,
    \item ``bursting a square'' is \cite[relation~2.10]{ckm}
      in the special case $r=s=l=1$.
  \end{itemize}
We will not reprint the full list of relations from \cite{ckm} here.
Their names for the relations not included in our subset are:
  \begin{itemize}
    \item two tag migrations \cite[relations~2.7 and 2.8]{ckm}
    \item square removal \cite[relation 2.9]{ckm}
    \item two digon removals
          \cite[relation 2.4 and 2.5]{ckm},
    \item the square switch \cite[relation 2.10]{ckm}.
  \end{itemize}
These relations can be derived from ours.  The proof assumes we can
divide by quantum integers, so it would not work in the important
setting where the scalars are $\mathbb{C}$ and $q \in \mathbb{C}$
is a root of unity.

\begin{prop}
\label{prop:equivalent}
  Our set of relations is equivalent to the full list of relations
  given in \cite{ckm}.
\end{prop}

\begin{proof}
  The first tag migration is $I=H$ in the special case $k+l+m=n$.
  The second follows from the first, using tag cancelation.

  Square removal in the case $s=1$ follows from from two applications
  of $I=H$, followed by bursting the resulting digon.  The general
  case then follows by induction on $s$.

  The first digon removal is square removal with two strands labeled
  $0$.  The second follows from the first by using tag cancellation
  and tag migration to change the orientations of some strands.

  To prove the square switch relation for $r=s=1$ and arbitrary
  $l>1$, consider the following web.
    \[
      \tikz {
	\draw[->-] (0,-5) -- (0,5); \draw[->-] (-30:5) arc (-30:30:5);
	\draw[-<-] (30:5) arc (30:90:5) node[pos=0.5,above]{$1$};
	\draw[->-] (90:5) arc (90:150:5) node[pos=0.5,above]{$1$};
	\draw[-<-] (150:5) arc (150:210:5); \draw[->-] (210:5) arc
	(210:270:5) node[pos=0.5,below]{$1$}; \draw[-<-] (270:5)
	arc (270:330:5) node[pos=0.5,below]{$1$}; \draw[->-] (0,-5)
	-- (0,5); \draw[->] (30:5) -- (30:10) node [right] {$l$};
	\draw[->] (150:5) -- (150:10) node [left] {$k$}; \draw[-<-]
	(210:5) -- (210:10) node [left] {$k$}; \draw[-<-] (330:5)
	-- (330:10) node [right] {$l$};
      }.
    \]
  This has two squares, either of which we can burst. Equating the
  two results, and using the quantum integer relation
    $$[k-1][l] - [l-1][k] = [k-l],$$
  gives the desired relation.

  It remains to prove that the general square switch relation follows
  from the special case $r=s=1$, together with square removal.  In
  the second remark of \cite{ckm}, they say ``there is an easy
  diagrammatic proof'' for this and other facts, ``or they can be
  proven as consequences of our main theorem''.  The diagrammatic
  proof requires a double induction on $r$ and $s$, and some messy
  calculations.  The less direct proof reduces to the fact that
  \cite[relation~4.1]{ckm} can be deduced from its special case
  $r=s=1$, together with \cite[relation~4.5]{ckm}, in their definition
  of $\dot{\mathcal{U}}_q(\mathfrak{gl}_2)$.  This is known, or can
  be derived by an algebraic calculation that mirrors the diagrammatic
  calculation.
\end{proof}

\section{Cobwebs}

Biologists classify spider webs into types.  If the diagrams in the
previous section are supposed to resemble the usual ``orb webs'', then
the diagrams introduced in this section more closely resemble ``cobwebs'',
since they have no vertices except for bivalent vertices and virtual
crossings.

Let $e_1,\dots,e_n$ be the standard basis vectors in $\mathbb{R}^n$.
The {\em roots} of the $SL_n$ root system are the vectors
  \[
    \alpha_{i,j} = e_i - e_j
  \]
for $i,j \in \{1,\dots,n\}$ with $i \neq j$.  Call $\alpha_{i,j}$ a
{\em positive root} if $i < j$.

\begin{figure}
\[
\begin{array}{OOO}
  \tikz {
    \draw[->-] (-4,-5) node [left] {$\alpha_{i,j}$} -- (0,0);
    \draw[->-] (4,-5) node [right] {$\alpha_{k,l}$} -- (0,0);
    \draw[->] (0,0) -- (-4,5) node [left] {$\alpha_{k,l}$};
    \draw[->] (0,0) -- (4,5) node [right] {$\alpha_{i,j}$};
  }
  & \tikz {
    \draw[->] (0,0) -- (0,5) node [right] {$\alpha_{i,j}$};
    \draw[->] (0,0) -- (0,-5) node [right] {$\alpha_{j,i}$};
    \draw (0,0) -- (2,0);
  }
  & \tikz {
    \draw[->-] (0,5) node [right] {$\alpha_{i,j}$} -- (0,0);
    \draw[->-] (0,-5) node [right] {$\alpha_{j,i}$} -- (0,0);
    \draw (0,0) -- (2,0);
  }
\end{array}
\]
\caption{Virtual crossings and bivalent vertices in a cobweb.}
\label{fig:cobgen}
\end{figure}

A {\em cobweb} is a collection of strands drawn in the plane.  Each
strand is oriented and is labeled by a root.  A cobweb can have
virtual crossings and bivalent vertices. A virtual crossing is a
point where a pair of strands pass through each other.  A bivalent
vertex has strands labeled $\alpha_{i,j}$ and $\alpha_{j,i}$,
oriented both in or both out, and has a tag that points to one side
of the vertex.  See Figure \ref{fig:cobgen}.

\begin{figure}
\[
\begin{array}{|T|T|}
  \hline
  {
    \phantom{ \tikz {\draw (0,-8) -- (0,8)}} 
    \tikz {
      \draw (0,-5) -- (0,5);
      \draw (0,0) -- (2,0);
    }
    = (-1) \tikz {
      \draw (0,-5) -- (0,5);
      \draw (0,0) -- (-2,0);
    } 
  } &
  {
    \tikz {
      \draw (0,-5) -- (0,5);
      \draw (0,-2) -- (2,-2);
      \draw (0,2) -- (2,2);
    }
    = \, \tikz {
      \draw (0,-5) -- (0,5);
    }
  } \\ \hline
  {
    \phantom{ \tikz {\draw (0,-8) -- (0,8)}} 
    \tikz {
      \draw[->] (5,0) arc (0:500:5)  node [left] {$\alpha_{i,j}$};
    }
    = q^{\sign(j - i)}
  } &
  {
    \tikz {
      \draw (-3,-5) .. controls (3,0) .. (-3,5);
      \draw (3,-5) .. controls (-3,0) .. (3,5);
    }
    = \tikz {
      \draw (-3,-5) .. controls (-2,-4) and (-2,4) .. (-3,5);
      \draw (3,-5) .. controls (2,-4) and (2,4) .. (3,5);
    }
  } \\ \hline
  {
    \phantom{ \tikz {\draw (0,-8) -- (0,8)}} 
    \tikz {
      \draw (-3,-5) -- (3,5);
      \draw (3,-5) -- (-3,5);
      \draw (-5,0) .. controls (0,3) .. (5,0);
    }
    = \tikz {
      \draw (-3,-5) -- (3,5);
      \draw (3,-5) -- (-3,5);
      \draw (-5,0) .. controls (0,-3) .. (5,0);
    }
  } &
  {
    \tikz {
      \draw (0,-5) -- (0,5);
      \draw (0,0) -- (2,0);
      \draw (-5,0) .. controls (0,3) .. (5,0);
    } =
    \tikz {
      \draw (0,-5) -- (0,5);
      \draw (0,0) -- (2,0);
      \draw (-5,0) .. controls (0,-3) .. (5,0);
    }
  } \\ \hline
  \multicolumn{2}{|F|}
  {
    \phantom{ \tikz {\draw (0,-8) -- (0,8)}} 
    \tikz {
      \draw[->] (-4,-5) -- (4,5) node [right] {$\alpha_{i,j}$};
      \draw[->] (4,-5) -- (-4,5) node [left] {$\alpha_{i,j}$};
    } =
    \tikz {
      \draw[->]
        (-4,-5) .. controls (0,0) .. (-4,5) node [left] {$\alpha_{i,j}$};
      \draw[->]
        (4,-5) .. controls (0,0) .. (4,5) node [right] {$\alpha_{i,j}$};
    }
  } \\ \hline
\end{array}
\]
\caption{Cobweb relations.
  Unlabeled strands can have any consistent orientations and labels.}
\label{fig:cobrel}
\end{figure}

The relations are shown in Figure \ref{fig:cobrel}.  These apply for
any way to fill in the omitted orientations and labels, where the
orientations and labels on a vertex or virtual crossing must be as
in the Figure~\ref{fig:cobgen}, and both sides of an equation must
have the same pattern of orientations and labels on their boundaries.
We call the relations:
  \begin{itemize}
    \item
      switching a tag,
    \item
      canceling a pair of tags,
    \item
      bursting a bubble,
    \item
      three detour moves: virtual Reidemeister moves two and three,
      and sliding a virtual crossing past a tag,
    \item
      smoothing a virtual crossing of two strands that have the same
      label.
  \end{itemize}

In general, diagrammatic relations can combine in complicated ways
to give unexpected relations, so we should check that our spider
of cobwebs is not trivial.

\begin{prop}
  \label{prop:nonzero}
  The spider of cobwebs is non-zero.
\end{prop}

\begin{proof}
  A \emph{closed} cobweb is a cobweb that has no endpoints on the
  boundary.  We will define a map that takes a closed cobweb to a
  monomial $\pm q^t$.  It is best here to think of each bivalent
  vertex as a place where a single strand changes its label and
  orientation.  With this convention, a closed cobweb consists of
  a collection of immersed loops, each having an even number of
  tags at which the orientation reverses and the label alternates
  between a root and its negative.

  Give each loop the overall orientation that is the same as all
  segments labeled by a positive root.  If a loop has no tags and
  is labeled by a negative root, give it the opposite orientation.
  With respect to these orientations, let $s$ be the number of tags
  that are on the left side of their loop, and let $t$ be the sum
  of the turning numbers of all loops.  Our desired monomial is
  then $(-1)^s q^t$.

  The resulting monomial is invariant under all cobweb relations.
  Thus we have a well-defined a surjective linear map from the
  vector space of closed cobwebs modulo cobweb relations onto the
  field of scalars.
\end{proof}

\begin{rmk}
  A similar argument shows that every cobweb is non-zero.  Furthermore,
  any two cobwebs with the same pattern of labels and orientations
  on the boundary are equal, up to multiplication by a scalar $\pm
  q^t$.  To prove this second fact, it helps to use the saddle move
  shown in Figure~\ref{fig:saddle}.
\end{rmk}

\section{Mapping webs to sums of cobwebs}

We will define a map that takes any given web to a sum of cobwebs.
We take a state sum, and describe how to build a cobweb corresponding
to each state.

A {\em state} on a web is a way to assign a subset of $\{1,\dots,n\}$
to each strand such that the following conditions are satisfied.
If a strand is labeled by the integer $k$, then it must be assigned
a set that has cardinality $k$.  The strands at a trivalent vertex
must be assigned sets $K$, $L$ and $K \cup L$, where $K$ and $L$
are disjoint.  The strands at a bivalent vertex must be assigned
$K$ and $L$ that are complementary subsets of $\{1,\dots,n\}$.

For every state, we define a corresponding cobweb as follows.  If a strand
is labeled by a set $K$ then replace it by a collection of parallel
strands, with one strand labeled $\alpha_{i,j}$ for every $i \in K$
and $j \not \in K$.  The orientation of these parallel strands should be
the same as that of the original strand.  The order they are arranged
in is not important.

Consider a trivalent vertex in which strands enter from the left
and right, and one strand exits to the top.  Suppose the strands
are assigned $K$ on the left, $L$ on the right, and $K \cup L$ on
the top, were $K$ and $L$ are disjoint subsets of $\{1,\dots,n\}$.
We replace this vertex with a cobweb that has the following strands:
  \begin{itemize}
  \item for every $i \in K$ and $j \not \in K \cup L$ there is
    a strand labeled $\alpha_{i,j}$ going from the left to the
    top,
  \item for every $i \in L$ and $j \not \in K \cup L$ there is
    a strand labeled $\alpha_{i,j}$ going from the right to
    the top, and
  \item for every $i \in K$ and $j \in L$ there are
    strands labeled $\alpha_{i,j}$ and $\alpha_{j,i}$ that come
    from the left and right to meet at a bivalent vertex with a tag
    pointing up.
  \end{itemize}
For a trivalent vertex with one strand oriented in and two strands
oriented out, simply reverse all of the arrows in the above
description.

Now consider a bivalent vertex with incoming strands entering from
the left and right, and a tag pointing up.  Suppose the strands are
assigned the sets $K$ on the left and $L$ on the right.  We map
this vertex to a cobweb with the following strands.  For every $i
\in K$ and $j \in L$, there are strands labeled $\alpha_{i,j}$ and
$\alpha_{j,i}$ coming from the left and right to meet at a bivalent
vertex with a tag pointing up.  For a bivalent vertex with two
strands oriented out, simply reverse all of the arrows in the above
description.

For any state on any web, the cobwebs at every strand and vertex
connect up to create a single cobweb.  The image of a web is the
sum over all possible states of the corresponding cobwebs.

\section{Mapping the $SL_n$ spider to the spider of cobwebs}

Let $\phi$ be the map defined in the previous section, extended
linearly to map linear combinations of webs to linear combinations
of cobwebs.  We will show this respects the defining relations on
webs, so is a well-defined map of spiders.  The proof of the following
theorem would work over any ring with an invertible element $q$.

\begin{thm}
\label{thm:welldef}
  $\phi$ is a well-defined map from the $SL_n$ spider
  to the spider of cobwebs.
\end{thm}

\begin{proof}
  For each web relation, replace each web with a state sum, then
  replace each state with a cobweb.  We must show that the resulting
  equation follows from the cobweb relations.

  Rather than calculate the entire state sum all at once, we will
  group terms according to their \emph{boundary data}.  That is,
  for each way to assign sets to the endpoints on the boundary of
  the disk, we sum over states that assign sets to all strands in
  a way that agrees with the given sets at the endpoints.

  For a given state, it is easiest to compute the resulting cobweb
  one pair of roots $\pm \alpha_{i,j}$ at a time.  We check cases
  based on which of the sets that label the edges contain $i$, $j$,
  both, or neither.  Cobweb strands with labels $\pm \alpha_{i,j}$
  appear wherever a web strand was assigned a set that contains
  exactly one of $i$ and $j$.

  \subsection{Strands labeled $0$ or $n$}

  A strand labeled $0$ must be assigned the empty set,
  which maps to an empty collection of cobweb strands.
  It can therefore be deleted.

  A strand labeled $n$ must be assigned the set $\{1,\dots,n\}$,
  which also maps to an empty collection of parallel cobweb strands.
  It can be deleted, leaving bivalent vertices that agree with the
  convention shown in Figure~\ref{fig:tag}.

  \subsection{Reversing all orientations}

  The cobweb relations are invariant under switching all orientations
  and mapping $q$ to $q^{-1}$.  The coefficients in the web relations
  are all invariant under mapping $q$ to $q^{-1}$, so all of our proofs
  will also apply when all orientations are reversed.

  \subsection{Tag operations}

  The tag switching and tag canceling relations for webs follow
  immediately from the same relations for cobwebs.  Note that a strand
  labeled $k$ is mapped to $k(n - k)$ parallel cobweb strands, so a tag
  switch in the web corresponds to $k(n - k)$ tag switches in the cobweb.

  \subsection{The $I = H$ relation}

  Assign the sets $K$, $L$, and $M$ to the strands labeled $k$, $l$
  and $m$, and $K \cup L \cup M$ to the strands at the top-right.  Here,
  $K$, $L$, and $M$ are disjoint and have cardinalities $k$, $l$, and $m$.
  We must then assign $L \cup M$ to the vertical strand in the $I$-shaped
  web, and $K \cup M$ to the horizontal strand in the $H$-shaped web.
  Each side of the equation is thus mapped to a single cobweb.

  Let $\alpha_{i,j}$ be a root.  If $i,j \in K$ then there is no
  strand labeled $\alpha_{i,j}$ or $\alpha_{j,i}$.  If $i \in K$
  and $j \in L$ then both sides of the equation have strands labeled
  $\alpha_{i,j}$ and $\alpha_{j,i}$ coming from the top-left and
  bottom-left, and meeting at a bivalent vertex with a tag to the
  right.  If $i \in K$ and $j \not\in K \cup L \cup M$ then both
  sides of the equation have a strand labeled $\alpha_{i,j}$ going
  from the top-left to the top-right. To be thorough, we could
  consider a total of sixteen cases like this, but they are all
  very similar.

  \subsection{Bursting a digon}

  Assign the set $K$ to the strands on the top and bottom of the
  webs.  On the left side of the equation, there are $k$ states
  consistent with this boundary data.  Namely, for any $a \in K$,
  assign $K \setminus \{a\}$ and $\{a\}$ to the strands on the left
  and right sides of the digon.

  The cobweb corresponding to this state has the following strands.
  For every $j \in K \setminus \{a\}$, there is a closed loop made
  by strands labeled $\alpha_{a,j}$ and $\alpha_{j,a}$ meeting at
  two bivalent vertices.  Every other strand goes from the bottom
  of the cobweb to the top.

  For every $j \in K \setminus \{a\}$, we can cancel a pair of tags
  to make a positively oriented bubble labeled $\alpha_{a,j}$.
  Bursting all of these bubbles gives $q$ to the power of the
  exponent:
    \[
      |\{j \in K \mid j > a\}| - |\{j \in K \mid j < a\}|.
    \]
  Now take the sum over all $a \in K$. The vertical strands of the
  cobwebs are the same in all terms, and the powers of $q$ sum to $[k]$.

  On the right side of the equation, there is only one state that
  is consistent with the boundary data.  Together with its coefficient
  $[k]$, we obtain the same as the left side of the equation.

  \subsection{Bursting a square}

  The possible types of boundary data fall into three cases.

  {\bf Case 1:}
    Let $K$ have cardinality $k$, and let $y \in K$.  Assign $K$ to the
    top-left and bottom-left strands of all webs, and $\{y\}$ to the
    top-right and bottom-right strands.

  On the left side of the equation, there are $k-1$ states consistent
  with the boundary data.  Namely, let $K' = K \setminus \{y\}$,
  and for any $a \in K'$, assign $\{a\}$ to the horizontal strands,
  $\{a,y\}$ to the right side of the square, and $K \setminus \{a\}$
  to the left.

  The cobweb corresponding to this state has the following
  strands.  For each $j \in K' \setminus \{a\}$, there is a closed
  loop made by strands labeled $\alpha_{a,j}$ and $\alpha_{j,a}$
  meeting at two bivalent vertices.  Every other strand goes
  from the bottom of the cobweb to the top, after canceling
  a pair of tags in the case of $\alpha_{y,a}$.

  For each $j \in K' \setminus \{a\}$, we can cancel a pair of tags
  to make a positively oriented bubble labeled $\alpha_{a,j}$.
  Bursting all of these bubbles gives $q$ to the power of the
  exponent
    \[
      |\{j \in K' \mid j > a\}| - |\{j \in K' \mid j < a\}|.
    \]
  Now take the sum over all $a \in K'$.  The vertical strands are the
  same for every term, and the powers of $q$ add up to $[k - 1]$.

  On the right side of the equation, the first ($I$-shaped) term
  has no state that is consistent with the boundary data, and the
  second has only one.  We get the coefficient $[k-1]$ times the
  cobweb with all strands going from the bottom to the top, which
  agrees with what we obtained for the left side of the equation.

  {\bf Case 2:}
    Let $K'$ have cardinality $k - 1$, and $x,y \not\in K'$.  Assign
    $K' \cup \{x\}$ and $\{y\}$ to the bottom-left and bottom-right
    strands of all webs. Assign $K' \cup \{y\}$ and $\{x\}$ to the
    top-left and top-right strands of all webs.

  On the left side of the equation, there is only one state that
  is consistent with the boundary data.  Namely, assign $\{y\}$ and
  $\{x\}$ to the top and bottom sides of the square, and assign
  $K'$ and $\{x,y\}$ to the left and right sides.  The corresponding
  cobweb has the following strands:
    \begin{itemize}
      \item for every $j \in K' \cup \{y\}$, there is a cup formed
        by strands labeled $\alpha_{x,j}$ and $\alpha_{j,x}$ that meet
        at a bivalent vertex with a downward tag,
      \item for every $j \in K' \cup \{x\}$, there is a cap formed
        by strands labeled $\alpha_{y,j}$ and $\alpha_{j,y}$ that meet
        at a bivalent vertex with an upward tag, and
      \item every other strand goes from the bottom of the cobweb to
        the top.
    \end{itemize}

  On the right side of the equation, the first term has one state
  that is consistent with the boundary data, and the second has
  none.  We get the same cobweb as we obtained for the left side
  of the equation.

  {\bf Case 3:}
    Let $K$ have cardinality $k$, and $y \not\in K$.  Assign $K$ to the
    top-left and bottom-left strands of all webs, and $\{y\}$ to the
    top-right and bottom-right strands.

  \begin{figure}
    \[
      \tikz {
        \draw[->-] (0,-5) node [left] {$\alpha_{x,y}$} ..
                   controls (4,-3) .. (5,-3);
        \draw (5,-3)--(5,-1.5);
        \draw[->-] (10,-5) node [right] {$\alpha_{y,x}$} ..
                   controls (6,-3) .. (5,-3);
        \draw[<-] (0,5) node [left] {$\alpha_{x,y}$} ..
                   controls (4,3) .. (5,3);
        \draw (5,3)--(5,1.5);
        \draw[<-] (10,5) node [right] {$\alpha_{y,x}$} ..
                   controls (6,3) .. (5,3);
      }
      =
      \tikz {
        \draw[->--] (0,-5) .. controls (4,3) .. (5,3);
        \draw (5,3)--(5,4.5);
        \draw[->--] (10,-5) .. controls (6,3) .. (5,3);
        \draw[<-] (0,5) ..  controls (4,-3) .. (5,-3);
        \draw (5,-3)--(5,-4.5);
        \draw[<-] (10,5) .. controls (6,-3) .. (5,-3);
      }
      =
      \tikz {
        \draw[->] (0,-5) .. controls (2,0) .. (0,5);
        \draw[->] (10,-5) .. controls (8,0) .. (10,5);
        \draw[->-] (5,-2) .. controls(3.9,-2) and  (3,-1.1) ..
                   (3,0) .. controls (3,1.1) and (3.9,2) .. (5,2);
        \draw[->-] (5,-2) .. controls(6.1,-2) and  (7,-1.1) ..
                   (7,0) .. controls (7,1.1) and (6.1,2) .. (5,2);
        \draw (5,2) -- (5,3.5);
        \draw (5,-2) -- (5,-3.5);
      }
      =
      q^{\sign(x-y)}
      \tikz {
        \draw[->] (0,-5) .. controls (2,0) .. (0,5);
        \draw[->] (10,-5) .. controls (8,0) .. (10,5);
      }
    \]
  \caption{A saddle move.}
  \label{fig:saddle}
  \end{figure}

  On the left side of the equation, there are $k$ states.  Namely,
  for any $a \in K$, assign $\{a\}$ to the top and bottom sides of
  the square, and assign $K \setminus \{a\}$ and $\{a,y\}$ to the
  left and right sides.  After canceling tags, the corresponding
  cobweb has the following strands:
    \begin{itemize}
      \item for $j \in K \setminus \{a\}$, there is a positively
        oriented bubble labeled $\alpha_{a,j}$,
      \item the strands labeled $\alpha_{a,y}$ and $\alpha_{y,a}$ form a
        cup and a cap, each having a bivalent vertex, and
      \item every other strand goes from the bottom of the
        cobweb to the top.
    \end{itemize}
  Use the \emph{saddle move} shown in Figure~\ref{fig:saddle} to
  replace the cup and cap with strands from the bottom of the cobweb
  to the top, with a coefficient $q^{\sign(a-y)}$.
  Bursting all bubbles now gives $q$ to the power of the exponent
    \[
      |\{j \in K \mid j > a\}| - |\{j \in K \mid j < a\}| + \sign(a - y).
    \]
  If $a$ is the smallest element of $K$
  then this exponent is $k-2$,
  or $k$ in the case $y<a$.
  As $a$ runs through the elements of $K$ in increasing order,
  the exponent decreases in steps of $2$,
  except that when and if $a$ jumps over $y$,
  there is a repeated exponent.
  In total, we get
    \[
      q^m + [k - 1],
      \text{ where }
      m = |\{j \in K \mid j > y\}| - |\{j \in K \mid j < y\}|.
    \]

  The first term on the right side of the equation has one state
  that is consistent with the boundary data.  The corresponding
  cobweb has the following strands.  For every $j \in K$, the strands
  labeled $\alpha_{j,y}$ and $\alpha_{y,j}$ form a cup and a cap
  with two bivalent vertices.  Every other strand goes from the
  bottom of the cobweb to the top.

  Use the saddle move again to make all strands go from the bottom
  of the cobweb to the top.  The total coefficient is $q^m$, with
  the same $m$ as above.

  The second term on the right side of the equation has one state
  that is consistent with the boundary data.  It has coefficient
  $[k-1]$, and the corresponding cobweb has all strands going from
  the bottom to the top.

  We see that both sides of the equation are equal to $q^m+[k-1]$
  times the same cobweb where all strands go from the bottom to the
  top.  This completes the proof.
\end{proof}

The fact that $\phi$ must be injective follows from the fact that
the $SL_n$ spider has no non-trivial proper monoidal ideals.  We
will explain this proof without going too deep into the category
theory.

\begin{thm}
\label{thm:injective}
  The map $\phi$ is injective.
\end{thm}

\begin{proof}
  For contradiction, suppose $v$ is a linear combination of webs
  that is non-zero in the $SL_n$ spider, but maps to zero in the
  spider of cobwebs.  We can assume $v$ is a linear combination of
  webs that have the same pattern of orientations and labels on
  their endpoints.

  Draw the webs that make up $v$ in a rectangle in such a way that
  all of their endpoints are on the top.  In the spider category,
  $v$ now represents a morphism whose domain is the unit object.

  Let $\Gamma_n$ be the functor defined in \cite{ckm} from the
  spider category to the representation category ${\mathcal{R}ep}(SL_n)$.
  By \cite[Theorem 3.3.1]{ckm}, $\Gamma_n$ is an equivalence of
  pivotal categories.  In particular, it is fully faithful, meaning
  it is a bijection on Hom spaces.

  Now, $\Gamma_n(v)$ is a morphism whose domain is the trivial
  one-dimensional representation.  Since $\Gamma_n$ is faithful,
  $\Gamma_n(v)$ is non-zero.  Since ${\mathcal{R}ep}(SL_n)$ is
  semisimple, some composition $f \circ \Gamma_n(v)$ is the identity
  morphism on the trivial representation.  Since $\Gamma_n$ is full,
  $f = \Gamma_n(u)$ for some $u$.  Since $\Gamma_n$ is faithful,
  $u \circ v$ is equal to the empty diagram.

  In the spider category, composition of webs is given by stacking,
  and this is extended bilinearly to linear combinations of webs.
  By assumption, $\phi(v)$ is zero in the spider of cobwebs.  Thus
  so is $\phi(u \circ v)$.  But $\phi(u \circ v)$ is the empty
  cobweb, so the empty cobweb is zero.  Every cobweb contains the
  empty cobweb, and hence is also zero, contradicting
  Proposition~\ref{prop:nonzero}.
\end{proof}

\section{Motivation}
\label{sec:motivation}

The point of a spider is to give a diagrammatic description of the
representation category of a quantum group.  One of the goals of
the so-called \emph{Kuperberg program} is to do as much as possible
using purely diagrammatic methods.  In particular, I would like to
construct spiders using as few prerequisites as possible.  Our
spider of cobwebs has the sort of simple construction I hoped for,
although unfortunately I do not know a good description of the image
of the $SL_n$ spider inside this simpler one.

We can at least partly explain the map from spiders to cobwebs in
terms of representation theory.  The functor $\Gamma_n$ takes a
strand labeled $k$ to the $k$th quantum exterior power of the
standard representation of $U_q(\mathfrak{sl}_n)$.  A state corresponds
to choosing specific vectors in each of these exterior powers.
Specifically, a strand labelled $K$ corresponds to the wedge of
basis vectors $e_i$ for $i \in K$, in increasing order.

For each set $K$ we can define a corresponding weight vector to be
$\sum_{i \in K} e_i$ projected onto the vector space spanned by the
roots, perpendicular to $(1,\dots,1)$.  For any state on a web, the
sum of weights corresponding to strands entering and leaving any
vertex the same.  Thus we could assign weights to the faces of the
web such that the weight assigned to an edge is the weight on the
face to its left minus the weight on the face to its right.

A \emph{graph planar algebra} is defined in \cite{jones} based on
this kind of state sum over ways of coloring the faces of a web.
We get an embedding the $SL_n$ spider in a version of the graph
planar algebra, where the graph is the weight lattice of $SL_n$.
We can evaluate such a state using \emph{Ocneanu cells}.  See
\cite{evanspugh} for the case $SL_3$, including generalizations at
special values of $q$.

Taking one more level of abstraction, we can specify a weight by
giving its inner product with every root.  In the cobweb, crossing
a strand from its right to its left corresponds to adding one to
the inner product with the root that labels that edge.  The cobweb
relations gives a way to evaluate a state on a closed cobweb, and
thus encode a simple version of Ocneanu cells.

Our map extends, with no extra effort, to webs with virtual crossings.
These webs satisfy detour moves, but not virtual Reidemeister one.
This leads to an invariant of ``rotational virtual'' knots and
links, which is presumably the same as the one defined by Kauffman
\cite{kauffman}.


\end{document}